\documentclass[12pt,a4paper]{article}
\usepackage{t1enc}
\usepackage[latin1]{inputenc}
\usepackage[english]{babel}
\usepackage{amsfonts}
\usepackage{amsmath}
\newcommand{\bp}{{\bf P}}
\newcommand{\be}{{\bf E}}

\def\begg{\begin{equation}}
\def\endd{\end{equation}}
\newtheorem{theorem}{Theorem}[section]
\newtheorem{corollary}{Corollary}[section]
\newtheorem{lemma}{Lemma}[section]

\newtheorem{remark}{Remark}[section]

\begin{document}
\centerline{\Large\bf On the behavior of random walk around heavy
points}

\bigskip\bigskip

\centerline{\it{ Dedicated to Professor Wolfgang Wertz on his
60-th birthday}}

\bigskip\bigskip
\renewcommand{\thefootnote}{1} \noindent \textbf{Endre
Cs\'{a}ki}\footnote{%
Research supported by the Hungarian National Foundation for Scientif\/ic
Research, Grant No. T 037886, T 043037 and K 061052.}\newline
Alfr\'ed R\'enyi Institute of Mathematics, Hungarian Academy of
Sciences,
Budapest, P.O.B. 127, H-1364, Hungary. E-mail address: csaki@renyi.hu

\bigskip

\renewcommand{\thefootnote}{2} \noindent \textbf{Ant\'{o}nia
F\"{o}ldes}%
\footnote{%
Research supported by a PSC CUNY Grant, No. 66494-0035.}\newline
Department of Mathematics, College of Staten Island, CUNY, 2800 Victory
Blvd., Staten Island, New York 10314, U.S.A. E-mail address:
foldes@mail.csi.cuny.edu

\bigskip

\noindent \textbf{P\'al R\'ev\'esz}$^1$ \newline
Institut f\"ur Statistik und Wahrscheinlichkeitstheorie, Technische
Universit\"at Wien, Wiedner Hauptstrasse 8-10/107 A-1040 Vienna,
Austria.
E-mail address: reveszp@renyi.hu

\bigskip \bigskip \bigskip

\noindent \textit{Abstract}: Consider a symmetric aperiodic random walk
in $Z^d$, $d\geq 3$. There are points (called heavy points) where the
number of visits by the random walk is close to its maximum. We
investigate the local times around these heavy points and show that they
converge to a deterministic limit as the number of steps tends to
infinity.

\bigskip


\bigskip


\bigskip
\noindent Running title: Random walk around heavy points

\renewcommand{\thesection}{\arabic{section}.}

\section{Introduction and main results}

\renewcommand{\thesection}{\arabic{section}} \setcounter{equation}{0} %
\setcounter{theorem}{0} \setcounter{lemma}{0}

Consider a random walk  $\{S_n\}_{n=1}^{\infty} $ starting at the origin
on the $d$-dimensional integer lattice $Z^d$, i.e.
$S_0=0$, $S_n=\sum_{k=1}^{n} X_k$, $n=1,2,\dots$, where
$X_k,\, k=1,2,\dots$ are i.i.d. random variables with distribution
\begg
\bp(X_1=x)=p(x),\quad x\in Z^d.
\label{px}
\endd

The random walk is called simple symmetric if $p(e_i)=1/(2d),\,
i=1,\ldots, 2d$, where $e_1,\ldots, e_d$ is a system of orthogonal unit
vectors in $Z^d$ and  $e_i=-e_{i-d},\, i=d+1,\ldots,2d$.

Denote by $Q$ the covariance matrix of $X_1$, and let $|Q|$ be its
determinant and let $Q^{-1}$ its inverse. Let
\begg
\Vert x\Vert^2:=xQ^{-1}x.
\label{norm}
\endd
For simple symmetric random walk $\Vert x\Vert^2
=|x|^2:=x_1^2+\cdots +x_d^2$, where
$x=(x_1,\ldots, x_d)$.

Recall the following definitions and basic properties from Spitzer \cite{Sp}.

A random walk is aperiodic if for
$$
R^+=\{x\in Z^d:\, \bp(S_n=x)>0\, \text{for some}\, n\geq 0\}
$$
we have
$$
\{x:\, x=y-z,\, \text{for some}\, y\in R^+,\, z\in R^+\}=Z^d.
$$
A random walk is strongly aperiodic if for each $x\in Z^d$ the
smallest subgroup  containing the set
$$
\{y:\, y=x+z,\, \text{where}\, p(z)>0\}
$$
is $Z^d.$  We assume throughout the paper that the random walk is
aperiodic (but not necessarily strongly aperiodic) and symmetric,
i.e. $p(x)=p(-x),\, x\in Z^d$.

For $d\geq 3$ the random walk is transient, i.e.
\begg
\gamma:=\bp(S_i\neq 0,\, i=1,2,\ldots)>0.
\label{gamma}
\endd

Define
\begg
\gamma_x:=\bp(S_i\neq x,\, i=1,2,\ldots),\quad x\in Z^d.
\label{gammax}
\endd

We shall impose the following moment conditions:
\begin{eqnarray}
\sum_{x\in Z^d} |x|^2p(x)&<&\infty, \qquad d=3,
\label{mom3} \\
\sum_{x\in Z^d} |x|^2\log (|x|+1)p(x)&<&\infty, \qquad d=4,
\label{mom4} \\
\sum_{x\in Z^d} |x|^{d-2}p(x)&<&\infty, \qquad d\geq 5,
\label{mom5}
\end{eqnarray}
where $|x|$ is the Euclidean distance.

The Green function is defined by
\begg
G(x):=\sum_{n=0}^\infty\bp(S_n=x),\quad x\in Z^d.
\label{green}
\endd

We have the identities
$$
\gamma=\frac1{G(0)},\qquad 1-\gamma_x=\frac{G(x)}{G(0)},\, x\neq 0.
$$

We need the following asymptotic property for the Green function in the
case of aperiodic random walk with mean 0, satisfying the moment
conditions (\ref{mom3}), (\ref{mom4}), (\ref{mom5}) for $d\geq 3$.

\begg
G(x)\sim c_d|Q|^{-1/2}\Vert x\Vert^{2-d},\quad |x|\to\infty
\label{greenas}
\endd
with some constant $c_d$. See Spitzer \cite{Sp}, p. 308 for $d=3$, p.
339, Problem 5 for $d>3$, or Uchiyama \cite{U} for strongly aperiodic
case and use Spitzer's trick (\cite{Sp}, p. 310) to reduce the aperiodic
case to strongly aperiodic case. For simple random walk see
R\'ev\'esz \cite{R}.

In this paper we are interested in studying local times of the random
walk defined by the number of visits as follows.
\begg
\xi(x,n):=\sum_{k=1}^{n}I\{S_k=x\},\quad
n=1,2,\ldots,\,
x\in Z^d,
\label{loct}
\endd
where $I\{A\}$ denotes the indicator of $A$.

Since the random walk is transient for $d\geq 3$, typically there is
only a finite number of visits to a fixed site, even for infinite time.
More precisely we have the distribution
\begg
\bp(\xi(0,\infty)=k)=\gamma(1-\gamma)^k, \quad k=0,1,2,\ldots
\label{loctd}
\endd
Cf. Erd\H os and Taylor \cite{ET60a} for simple random walk. The general
case is similar.

There are however (random) points where the random walk
accumulates a higher number of visits. Consider the maximal local
time \begg \xi(n):=\max_{x\in Z^d}\xi(x,n), \quad n=1,2,\ldots
\label{maxxi}
\endd
and also
\begg
\eta(n):=\max_{0\leq j\leq n}\xi(S_j,\infty), \quad n=1,2,\ldots
\label{maxeta}
\endd

Erd\H os and Taylor \cite{ET60a} proved for simple random walk and
$d\geq 3$
\begg
\lim_{n\to\infty}\frac{\xi(n)}{\log n}=\lambda:
=-\frac1{\log(1-\gamma)}\hspace{1cm} \mathrm{a.s.}
\label{limxi}
\endd

\bigskip Following the proof of Erd\H{o}s and Taylor, without any new
idea, one can prove that (\ref{limxi}) holds for general aperiodic
random walk and also
\begin{equation}
\lim_{n\to\infty}\frac{\eta(n)}{\log n}=\lambda \hspace{1cm}
\mathrm{a.s.}
\label{lb}
\end{equation}

For general treatment of similar strong theorems for local and
occupation times see \cite{CsFRRS}.

(\ref{limxi}) means that there are sites where the local time up
to time $n$ is around $\lambda\log n$. These will be called heavy
points. We are interested in the problem what happens around these
heavy points. We may ask whether it is possible that in a close
neighborhood of a heavy point there is another heavy point? Or an
empty point (not visited at all up to time $n$)? We shall see that
the answers for both questions happen to be negative.

In \cite{CsFR06} we investigated the joint asymptotic behavior of local
times of two neighboring sites for simple random walk and found that the
vector
$$
\left(\frac{\xi(x,n)}{\log n},\, \frac{\xi(x+e_1,n)}{\log n}\right)
$$
is essentially in the domain
$$
\{y\geq 0, z\geq 0:\, -(y+z)\log (y+z)+y\log y+z\log z-(y+z)\log\alpha
\leq 1\},
$$
where
$$
\alpha:=\frac{1-\gamma}{2-\gamma}.
$$

One can see that the only point in this domain with $y=\lambda$ is
$z=\lambda(1-\gamma)$, which tells us that if a point is heavy, i.e. its
local time is around $\lambda\log n$, then the local time of any of its
neighbors should be around $\lambda(1-\gamma)\log n$, i.e. cannot
fluctuate too much, at least asymptotically. We say that the local time
around a heavy point is asymptotically deterministic. Our concern is to
investigate this phenomenon further and determine the asymptotic value
of local times of sites $x$ with $\Vert x\Vert\leq r_n$, where $r_n$ may
tend to infinity at a certain rate.

Define
\begg
m_x=\left\{
\begin{array}{ll}
& 1\quad\quad\quad\ {\rm if}\quad x=0, \\
& \frac{(1-\gamma_x)^2}{1-\gamma}\quad {\rm if}\quad x\neq 0.
\end{array}
\right.
\label{mx}
\endd
$m_x$ is, in fact, the expectation of the local time at $x$ between two
consecutive returns to zero (see Remark 2.1).

We shall consider the "balls" (which are, in fact, ellipsoids in
Euclidean space)
\begg
B(r)=\left\{x:\, \Vert x\Vert
\leq r\right\},
\label{ball}
\endd
where $\Vert x\Vert$ is defined by (\ref{norm}).

\begin{theorem} Let $d\geq 5$ and $k_n=(1-\delta_n)\lambda\log n$. Let
$r_n>0$ and $\delta_n>0$ be selected such that  $\delta_n$ is
non-increasing, $r_n$ is non-decreasing, and for any $c>0,$ let
$r_{[cn]}/r_n<C$ with some $C>0$ and for \begg
\beta_n:=r_n^{2d-4}\frac{\log\log n}{\log n} \label{beta1}
\endd
\begg
\lim_{n\to\infty}\beta_n=0,\qquad \lim_{n\to\infty}\delta_n
r_n^{2d-4}=0.
\endd
Define the random set of points
\begg
{\cal A}_n=\{z\in Z^d: \xi(z,n)\geq k_n\}.
\endd
Then we have for symmetric aperiodic random walk
\begg
\lim_{n\to\infty}\sup_{z\in {\cal A}_n}\sup_{x\in B(r_n)}
\left|\frac{\xi(z+x,n)}{m_x\lambda\log n}-1\right|=0\quad{\rm a.s.}
\endd
\end{theorem}

\begin{theorem} Let $d\geq 3$ and $k_n=(1-\delta_n)\lambda\log n$. Let
$r_n>0$ and $\delta_n>0$ be selected such that $\delta_n$ is
non-increasing, $r_n$ is non-decreasing, and for any $c>0,$ let
$r_{[cn]}/r_n<C$ for some $C>0$ and for \begg
\beta_n:=r_n^{2d-4}\frac{\log\log n}{\log n} \label{beta2}
\endd
\begg
\lim_{n\to\infty}\beta_n=0,\qquad \lim_{n\to\infty}\delta_n
r_n^{2d-4}=0.
\endd
Define the random set of indices
\begg
{\cal B}_n=\{j\leq n: \xi(S_j,\infty)\geq k_n\}.
\endd
Then we have for symmetric aperiodic random walk
\begg
\lim_{n\to\infty}\sup_{j\in {\cal B}_n}\sup_{x\in B(r_n)}
\left|\frac{\xi(S_j+x,\infty)}{m_x\lambda\log n}-1\right|=0\quad{\rm
a.s.}
\endd
\end{theorem}

\begin{remark} For a given $\omega$, ${\cal A}_n$ or ${\cal B}_n$ can
be empty. In this case $\sup_{z\in {\cal A}_n}$ or $\sup_{j\in {\cal
B}_n}$ is automatically considered to be $0$.
\end{remark}

\begin{corollary}
Let $A\subset Z^d$ be a fixed set.

{\rm (i)}
If $d\geq 5$ and $z_n\in {\cal A}_n$, then
$$
\lim_{n\to\infty}\frac{\sum_{x\in A}\xi(x+z_n,n)}{\log n}
=\lambda\sum_{x\in A}m_x \quad {\rm a.s.}
$$

{\rm (ii)}
If $d\geq 3$ and $j_n\in {\cal B}_n$, then
$$
\lim_{n\to\infty}\frac{\sum_{x\in A}\xi(x+S_{j_n},\infty)}{\log n}
=\lambda\sum_{x\in A}m_x \quad {\rm a.s.}
$$
\end{corollary}

\bigskip
 From our Theorems it is obvious that the critical case is around $r_n\sim
(\log n)^{1/(2d-4)}$. It follows that for smaller $r_n$ the ball
$S_j+B(r_n)$ is completely covered for $j\in {\cal B}_n$ with
probability 1. We have the following Corollary.

\begin{corollary}
For $j\in {\cal B}_n$ let $R(n,j)$ denote the largest number such that
$S_j+B(R(n,j))$ is completely covered by the random walk
$S_0,S_1,S_2,\ldots$, i.e. $\xi(S_j+x,\infty)>0$, $x\in B(R(n,j))$. Then
for any $\varepsilon>0$ we have $R(n,j)\geq (\log
n)^{(1-\varepsilon)/(2d-4)}$ almost surely.
\end{corollary}

We conjecture that for $j\in{\cal B}_n$ we have
$R(n,j)\leq (\log n)^{(1+\varepsilon)/(2d-4)}$. Our next result is one
step in this direction, showing that in Theorems 1.2 the power $1/(2d-4)$
of $\log n$ cannot be improved in general.

\begin{theorem}
For simple symmetric random walk let $\{x_n\}$ be a sequence such that
$|x_n|\sim c(\log n)^{1/(2d-4)}$ for some $c>0$. Then with probability
one there exist infinitely many $n$ such that
$$
\xi(S_n,\infty)\geq \lambda\left(\log n
+\left(\frac{d-4}{d-2}-\varepsilon\right)\log\log n\right),\quad
\xi(S_n+x_n,\infty)=0.
$$
Consequently, $n\in{\cal B}_n$ and $R(n,n)\leq c(\log n)^{1/(2d-4)}$
infinitely often with probability one.
\end{theorem}

\renewcommand{\thesection}{\arabic{section}.}

\section{Preliminary facts and results}

\renewcommand{\thesection}{\arabic{section}} \setcounter{equation}{0}
\setcounter{theorem}{0} \setcounter{lemma}{0}

First we present some more notations. For $x\in Z^d$
let $T_x$ be the first hitting time of the point $x$, i.e.
$T_x=\min\{i\geq 1:S_i=x\}$ with the convention that
$T_x=\infty$ if there is no $i$ with $S_i=x$. Denote $T_0=T$.

Introduce further
\begin{eqnarray}
q_x&:=&\bp(T<T_x),\label{aq} \\
s_x&:=&\bp(T_x<T).\label{as}
\end{eqnarray}
In words, $q_x$ is the probability that the random walk, starting from
$0$, returns to $0$, before hitting $x$ (including $T<T_x=\infty$), and
$s_x$ is the probability that the random walk, starting from $0$, hits
$x$, before returning to $0$ (including $T_x<T=\infty$).

Now we give the joint distribution of $\xi(0,\infty)$ and
$\xi(x,\infty)$ in the following form.

\begin{lemma} For $x\neq 0$, $v<\log(1/(1-\gamma))$, $k=0,1,2,\ldots$
\begg
\be(e^{v\xi(x,\infty)};\, \xi(0,\infty)=k)=
\left(q_x+\frac{s_x^2e^v}{1-q_xe^v}\right)^k(1-q_x-s_x)
\left(1+\frac{s_xe^v}{1-q_xe^v}\right)
\label{momgen}
\endd
\begg
=\gamma(1-\gamma)^k\left(\varphi(v)\right)^k\psi(v),
\label{fipsi}
\endd
where
\begg
\varphi(v):=\frac{1-\frac{(1-\gamma)^2-(1-\gamma_x)^2}
{\gamma(1-\gamma)}(e^v-1)}
{1-\frac{1-\gamma-(1-\gamma_x)^2}{\gamma}(e^v-1)},
\endd
\begg \psi(v):=\frac{1-\frac{\gamma_x-\gamma} {\gamma}(e^v-1)}
{1-\frac{1-\gamma-(1-\gamma_x)^2}{\gamma}(e^v-1)}.
\endd
\end{lemma}

\noindent{\bf Proof.} Observe that
\begg
\bp \left(\sum_{n=1}^T I\{S_n=x\}=j,T<\infty\right)=\left\{
\begin{array}{ll}
q_x& \quad\mathrm{if \, \,}j=0 , \\
s_x^2q_x^{j-1} & \quad\mathrm{if \, \,} j=1,2,...
\end{array}
\right.
\label{finite}
\endd
and
\begg
\bp \left(\sum_{n=1}^T I\{S_n=x\}=j,\, T=\infty\right)=\left\{
\begin{array}{ll}
1-q_x-s_x& \quad\mathrm{if \, \,}j=0 , \\
s_x(1-q_x-s_x)q_x^{j-1} & \quad\mathrm{if \, \,} j=1,2,...
\end{array}
\right.
\label{infinite}
\endd

Obviously
$$
\xi(x,\infty)=Z_1+\ldots+Z_{\xi(0,\infty)}+\hat Z,
$$
where $Z_1,\ldots, Z_{\xi(0,\infty)}$ are the local times of $x$
between consecutive returns to $0$ and $\hat Z$ is the local time
of $x$ after the last return to zero. Hence (\ref{momgen}) follows
from (\ref{finite}) and (\ref{infinite}). (\ref{fipsi}) can be
obtained by using
\begin{eqnarray}
q_x&=&1-\frac{{\gamma}}{1-(1-{\gamma}_x)^2},  \label{qu} \\
s_x&=&(1-{\gamma}_x)(1-q_x).  \label{es}
\end{eqnarray}
(Cf. \cite{CsFR05} or \cite{R} for simple random walk, the
general case being similar).

\begin{remark} It is easy to  see that our condition $v<\log(1/(1-\gamma))$
implies $q_xe^v<1$, needed to obtain {\rm (\ref{momgen})}.
Furthermore
$$
\varphi(v)=\be\left(e^{v\sum_{n=1}^T I\{S_n=x\}}\mid T<\infty\right),
$$
$$
\psi(v)=\be\left(e^{v\sum_{n=1}^T I\{S_n=x\}}\mid T=\infty\right)
$$
and
$$
m_x=\be\left(\sum_{n=1}^T I\{S_n=x\}\mid T<\infty\right).
$$
\end{remark}

Further properties of $q_x$ and $s_x$ for simple symmetric random walk
is given in the next Lemma.

\begin{lemma}
For simple symmetric random walk and $x\in Z^d$
\begin{eqnarray}
\gamma_x&\geq& \gamma, \label{xe} \\
\frac{1-\gamma}{2-\gamma}&\leq& q_x\leq 1-\gamma, \label{qe} \\
1-q_x-s_x&\geq& \frac{\gamma}{2-\gamma},\label{qse}\\
q_x(n)&:=&\bp(T<\min(n,T_x))=q_x+\frac{O(1)}{n^{d/2-1}}.\label{qxn}
\end{eqnarray}
\end{lemma}

\noindent {\bf Proof.}
For (\ref{xe}) see \cite{CsFR05}, Lemma 2.4 and for (\ref{qxn}) see
\cite{CsFR05}, Lemma 2.5. (\ref{qe}) and (\ref{qse}) can be easily
obtained from (\ref{qu}), (\ref{es}) and (\ref{xe}).

The next result gives an estimation of $\varphi$ and $\psi$, where the
error term is uniform in $x$.

\begin{lemma} For $\log(1-\gamma(1-\gamma))<v<\log(1+\gamma(1-\gamma))$
we have
\begg
\varphi(v)=\exp(m_x (v+O(v^2))),\quad v\to 0,
\label{fi}
\endd
where $O$ is uniform in $x$,
\begg
\psi(v)\leq \frac{1+|e^v-1|}{1-|e^v-1|/\gamma}.
\label{psiest}
\endd
\end{lemma}

\noindent{\bf Proof.} Write
$$
\varphi(v)=\frac{1-u}{1-y}
$$
with
$$
u=\frac{(1-\gamma)^2-(1-\gamma_x)^2}{\gamma(1-\gamma)}(e^v-1),
\quad y=\frac{1-\gamma-(1-\gamma_x)^2}{\gamma}(e^v-1).
$$
Then it is easy to see that
$$
y-u=m_x(e^v-1),
$$
and
$$
|u|\leq\frac{|e^v-1|}{\gamma(1-\gamma)}, \quad
|y|\leq\frac{|e^v-1|}{\gamma(1-\gamma)}.
$$

By Taylor series
$$
\log\frac{1-u}{1-y}=\log(1-u)-\log(1-y)=y-u+\frac{y^2-u^2}2
+\frac{y^3-u^3}3+\ldots
$$
$$
=(y-u)\left(1+\frac{y+u}2+\frac{y^2+uy+u^2}3+\ldots\right).
$$

Since $e^v-1=v+O(v^2)$, we have
$$
\left| \log\frac{1-u}{1-y}-m_x(e^v-1)\right|
\leq m_x|e^v-1|
\left(\frac{|e^v-1|}{\gamma(1-\gamma)}+
\left(\frac{|e^v-1|}{\gamma(1-\gamma)}\right)^2+\ldots\right)=m_xO(v^2),
$$
where $O$ is independent of $x$. Hence (\ref{fi}) follows.
(\ref{psiest}) is obvious.

\renewcommand{\thesection}{\arabic{section}.}

\section{Proof of Theorem 1.2}

\renewcommand{\thesection}{\arabic{section}} \setcounter{equation}{0}
\setcounter{theorem}{0} \setcounter{lemma}{0}

Observe that $k_n\sim \lambda\log n$. Let $n_\ell=[e^\ell]$,
 and define the events
$$
A_j=\left\{\xi(S_j,\infty)\geq k_{n_\ell},\, \sup_{x\in
B(r_{n_{\ell+1}})}\left(\frac{\xi(S_j+x,\infty)}
{m_xk_{n_\ell}}-1\right)\geq\varepsilon\right\}
$$
$$
\bp\left(\bigcup_{j=0}^{n_{\ell+1}}A_j\right)
\leq\sum_{j=0}^{n_{\ell+1}}\bp(A_j)
\leq \sum_{j=0}^{n_{\ell+1}}\sum_{x\in
B(r_{n_{\ell+1}})}\bp(A_j^{(x)}),
$$
where
$$
A_j^{(x)}=\left\{\xi(S_j,\infty)\geq k_{n_\ell},\,
\xi(S_j+x,\infty)\geq (1+\varepsilon)m_x k_{n_\ell}\right\}.
$$

Consider the random walk obtained by reversing the original walk at $S_j$,
i.e. let $S_i':=S_{j-i}-S_j$, $i=0,1,\ldots,j$ and extend it to infinite
time, and also the forward random walk $S_i'':=S_{j+i}-S_j$,
$i=0,1,2,\ldots$ Then $\{S_0',S_1',\ldots\}$ and
$\{S_0'',S_1'',\ldots\}$ are independent random walks and so are their
respective local times $\xi'$ and $\xi"$. Moreover,

$$
\xi(S_j,\infty)=\xi"(0,\infty)+\xi(S_j,j)\leq
\xi"(0,\infty)+\xi'(0,\infty)+1,
$$
$$
\xi(S_j+x,\infty)=\xi"(x,\infty)+\xi(S_j+x,j)\leq
\xi"(x,\infty)+\xi'(x,\infty).
$$
Here $\xi'$ and $\xi"$ are independent and have the same distribution as
$\xi$.

Hence
\begin{eqnarray*}
\bp(A_j^{(x)}) &\leq&\bp(\xi"(0,\infty)+\xi'(0,\infty)\geq
k_{n_\ell}-1,\, \xi"(x,\infty)+\xi'(x,\infty)\geq
(1+\varepsilon)m_x k_{n_\ell})\\
&=&\sum \bp(\xi"(0,\infty)=k_1,\xi'(0,\infty)=k_2,
\xi"(x,\infty)+\xi'(x,\infty)\geq (1+\varepsilon) m_x k_{n_\ell}),
\end{eqnarray*}
where the summation goes for $k_1+k_2\geq k_{n_\ell}-1$. Using
exponential Markov inequality, Lemma 2.1, independence of $\xi"$ and
$\xi'$ and elementary calculus, we get
\begin{eqnarray*}
\bp(A_j^{(x)}) &\leq& \sum \be
\left(e^{v(\xi"(x,\infty)+\xi'(x,\infty))},\xi"(0,\infty)=k_1,
\xi'(0,\infty)=k_2\right)
e^{-v(1+\varepsilon)m_x k_{n_\ell}} \\
&=&\sum (\varphi(v))^{k_1+k_2}\gamma^2(1-\gamma)^{k_1+k_2}\psi^2(v)
e^{-v(1+\varepsilon)m_xk_{n_\ell}} \\
&=&\gamma^2\psi^2(v)e^{-v(1+\varepsilon)m_xk_{n_\ell}}\sum
(\varphi(v)(1-\gamma))^{k_1+k_2} \\
&=&\gamma^2\psi^2(v)e^{-v(1+\varepsilon)m_xk_{n_\ell}}
(\varphi(v)(1-\gamma))^{k_{n_\ell}} \\
&\times&\left(\frac{k_{n_\ell}}{\varphi(v)(1-\gamma)(1-\varphi(v)(1-\gamma))}
+\frac1{(1-\varphi(v)(1-\gamma))^2}\right).
\end{eqnarray*}

By (\ref{fi}) we obtain for all $j$
\begin{eqnarray*}
\bp(A_j^{(x)})&\leq& \gamma^2\psi^2(v)
\left(\frac{k_{n_\ell}}{\varphi(v)(1-\gamma)(1-\varphi(v)(1-\gamma))}
+\frac1{(1-\varphi(v)(1-\gamma))^2}\right) \\
&\times& e^{-m_xvk_{n_\ell}(\varepsilon+O(v))}(1-\gamma)^{k_{n_\ell}}.
\end{eqnarray*}
Choose $v_0>0$ small enough such that
$$
\varepsilon+O(v_0)>0,\quad e^{v_0}<1+\gamma(1-\gamma),\quad
\varphi(v_0)<\frac1{1-\gamma}.
$$
Using $x\in B(r_{n_{\ell+1}})$ and (\ref{greenas}) we get
$$
m_xk_{n_\ell}=\frac{(1-\gamma_x)^2}{1-\gamma}
(\lambda\log n_\ell(1-\delta_{n_\ell}))\geq
\frac{C_1(1-\delta_{n_\ell})\log n_\ell}{\Vert x\Vert^{2d-4}}
\geq\frac{C_1(1-\delta_{n_\ell})\log n_\ell}{r_{n_{\ell+1}}^{2d-4}},
$$
where here and in the sequel $C_1,C_2,\ldots$ will denote positive
constants whose values are unimportant in our proofs.

By the above assumptions
\begin{eqnarray*}
\bp(A_j^{(x)})&\leq& C_2k_{n_\ell}
e^{-m_xv_0 k_{n_\ell}(\varepsilon+O(v_0))}(1-\gamma)^{k_{n_\ell}} \\
&\leq& C_2k_{n_\ell}
\exp\left(-(1-\delta_{n_\ell})\log n_\ell
\left(\frac{C_3}{r_{n_{\ell+1}}^{2d-4}}+1\right)\right).
\end{eqnarray*}

Hence
\begin{eqnarray*}
\sum_{j=0}^{n_{\ell+1}}\sum_{x\in
B(r_{n_{\ell+1}})}\bp(A_j^{(x)})
&\leq& C_4n_{\ell+1}r_{n_{\ell+1}}^d
k_{n_\ell}
\exp\left(-(1-\delta_{n_\ell})\log n_\ell
\left(\frac{C_3}{r_{n_{\ell+1}}^{2d-4}}+1\right)\right) \\
&\leq& C_4\frac{n_{\ell+1}}{n_\ell}k_{n_\ell}r_{n_{\ell+1}}^d
\exp\left(-\frac{C_3\log n_\ell}{r_{n_{\ell+1}}^{2d-4}}+
\delta_{n_\ell}\log n_\ell\right) \\
&=& C_4\frac{n_{\ell+1}}{n_\ell}k_{n_\ell}r_{n_{\ell+1}}^d
\exp\left(-\frac{\log n_\ell}{r_{n_{\ell}}^{2d-4}}
\left(C_3\left(\frac{r_{n_\ell}}{r_{n_{\ell+1}}}\right)^{2d-4}
-\delta_{n_\ell}r_{n_\ell}^{2d-4}\right)\right) \\
&\leq& C_4\frac{n_{\ell+1}}{n_\ell}k_{n_\ell}r_{n_{\ell+1}}^d
\exp\left(-C_5\frac{\log n_\ell}{r_{n_{\ell}}^{2d-4}}\right)
\leq C_6(\log n_\ell)^{3-\frac{C_7}{\beta_{n_\ell}}},
\end{eqnarray*}
where in the last two lines we used the conditions of the Theorem for
$r_n$ and $\delta_n$. Consequently
$$
\bp(\bigcup_{j=0}^{n_{\ell+1}} A_j)\leq
\sum_{j=0}^{n_{\ell+1}}\sum_{x\in
B(r_{n_{\ell+1}})}\bp(A_j^{(x)})
\leq C_6 \ell^{3-\frac{C_7}{\beta_{n_\ell}}}\leq \frac{C_6}{\ell^2}
$$
for large enough $\ell$ which is summable in $\ell$. By
Borel-Cantelli lemma for large $\ell$ if $\xi(S_j,\infty)\geq
k_{n_\ell}$, then $\xi(S_j+x,\infty)\leq (1+\varepsilon)m_x
k_{n_\ell}$ for all $x\in B(r_{n_{\ell+1}})$.

Let now $n_\ell\leq n<n_{\ell+1}$ and $x\in B(r_{n_{\ell+1}})$.
$\xi(S_j,\infty)\geq k_n, j\leq n$ implies $\xi(S_j,\infty)\geq
k_{n_\ell}$, i.e.
\begin{equation}
\xi(S_j+x,\infty)\leq (1+\varepsilon)m_xk_{n_\ell}\leq
(1+\varepsilon)m_xk_n. \label{xisx}
\end{equation}

The lower bound is similar, with slight modifications. We call $S_j$ new
if $S_i\neq S_j,\, i=1,2,\ldots, j-1$. Define the events
$$
D_j=\left\{\xi(S_j,\infty)\geq k_{n_\ell},\, \sup_{x\in
B(r_{n_{\ell+1}})}\left(1-\frac{\xi(S_j+x,\infty)}
{m_xk_{n_{\ell+1}}}\right)\geq\varepsilon\right\},
$$
$$
D_j^{(x)}=\{S_j\, {\rm new},\, \xi(S_j,\infty)\geq k_{n_\ell},
\xi(S_j+x,\infty)\leq (1-\varepsilon)m_xk_{n_{\ell+1}}\}.
$$
Observe that
$$
\bigcup_{\{j:0\leq j\leq n_{\ell+1}\}}D_j=
\bigcup_{\{j:0\leq j\leq n_{\ell+1},\, S_j\, {\rm new}\}}D_j.
$$
Considering again the forward random walk, we have
$$
\xi(S_j,\infty)=\xi"(0,\infty)+1,\,\,\, \xi(S_j+x,\infty)\geq
\xi"(x,\infty).
$$
Hence by Markov's inequality
\begin{eqnarray*}
\bp(D_j^{(x)})&\leq& \sum_{k=k_{n_\ell}-1}^\infty
\bp(\xi"(0,\infty)=k, \xi"(x,\infty)\leq (1-\varepsilon)m_xk_{n_{\ell+1}}) \\
&\leq& \sum_{k=k_{n_\ell}-1}^\infty (\varphi(-v)(1-\gamma))^k
\psi(-v)\exp(v(1-\varepsilon)m_xk_{n_{\ell+1}}) \\
&\leq&\frac{\psi(-v)}{(1-\gamma)\varphi(-v)(1-(1-\gamma)\varphi(-v))}
((1-\gamma)\varphi(-v))^{k_{n_\ell}}e^{v(1-\varepsilon)m_xk_{n_{\ell+1}}}.
\end{eqnarray*}

Proceeding as above we finally conclude after somewhat simpler
calculations than the previous one, that for large enough $n$,
$\xi(S_j,\infty)\geq k_n$ implies $\xi(S_j+x,\infty)\geq
(1-\varepsilon)m_xk_n$.

This, combined with (\ref{xisx}) completes the proof of Theorem 1.2.

\renewcommand{\thesection}{\arabic{section}.}

\section{Proof of Theorem 1.1}

\renewcommand{\thesection}{\arabic{section}} \setcounter{equation}{0}
\setcounter{theorem}{0} \setcounter{lemma}{0}

\begin{lemma}
Let $d\geq 5$, $\frac{2}{d-2}<\alpha<1$, $j\leq n-n^\alpha$,
$|x|\leq \log n$. Then with probability 1 there exists an
$n_0(\omega)$ such that for $n\geq n_0$ we have
$$
\xi(S_j+x,n)=\xi(S_j+x,\infty).
$$
\end{lemma}

\noindent
{\bf Proof.} The proof is essentially the same as that of Theorem 1
(iii) in Erd\H os and Taylor \cite{ET60b}.

Let
$$
n_{k+1}=n_k+\left[\frac12 n_k^\alpha\right].
$$
$$
A_k=\bigcup_{j\leq n_k} \quad\bigcup_{\ell\geq n_k+[\frac12
n_{k-1}^\alpha]} \quad\bigcup_{x\in B(\log (2
n_{k+1}))}\{S_\ell-S_j=x\}.
$$
For aperiodic random walk we have (cf. Jain and Pruitt \cite{JP})
\begin{equation}
\bp(S_n=x)\leq C_8n^{-d/2}
\label{jp}
\end{equation}
for all $x\in Z^d$ and $n\geq 1$ with some constant $C_8$.

Using the fact that $B(\log (2n_{k+1}))$ contains less than
$C_9(\log n_{k+1})^d$ points,
\begin{eqnarray}
\bp(A_k)&\leq&
C_9(\log n_{k+1})^d \sum_{j=0}^{n_k}
\sum_{\ell=n_k+[\frac12n_{k-1}^\alpha]}^\infty
\frac{C_8}{(\ell-j)^{d/2}} \nonumber \\
&\leq& \sum_{j=0}^{n_k} \frac{C_{10}(\log
n_{k+1})^d}{(n_k+[\frac12n_{k-1}^\alpha]-j)^{d/2-1}} \leq \frac
{C_{10}(\log n_{k+1})^d}{n_{k-1}^{\alpha(d/2-2)}} \leq \frac
{C_{11}(\log n_{k-1})^d}{n_{k-1}^{\alpha(d-4)/2}}. \label{pak1}
\end{eqnarray}
We will show now that $\sum_k\bp(A_k)$ converges.
\begin{eqnarray}
\sum_{n=1}^\infty \frac{(\log n)^d}{n^{\alpha(d-2)/2}}&\geq&
\sum_k\sum_{n=n_k+1}^{n_{k+1}}\frac{(\log n)^d}{n^{\alpha(d-2)/2}}
\geq C_{12}\sum_k \frac{n_{k+1}-n_k}{n_{k+1}^{\alpha(d-2)/2}}(\log
n_{k+1})^d \nonumber \\
&\geq& C_{12}\sum_k \frac{\frac12
n_k^\alpha}{n_{k+1}^{\alpha(d-2)/2}} (\log n_{k+1})^d =
C_{13}\sum_k \frac{(\log n_{k+1})^d}{n_{k+1}^{\alpha(d-4)/2}}
\left(\frac{n_k}{n_{k+1}}\right)^\alpha. \label{pak2}
\end{eqnarray}
Observe that
$$
\left(\frac{n_k}{n_{k+1}}\right)^\alpha=
\left(\frac{n_k}{n_k+[\frac12n_k^\alpha]}\right)^\alpha\to 1,\quad
k\to\infty.
$$
Since
$$
\sum_{n=1}^\infty \frac{(\log n)^d}{n^{\alpha(d-2)/2}}
$$
converges, (\ref{pak1}) and (\ref{pak2}) imply the convergence of
$\sum_k\bp(A_k)$. By Borel-Cantelli lemma, if $k$ is big enough,
the tube of radius $\log ( 2 n_{k+1})$ around the path $\{S_j,\,
j=1,2,\ldots, n_k\}$ is disjoint from the path $\{S_\ell,\,
\ell=n_k+[\frac12n_{k-1}^\alpha],\ldots\}$.

To finish the proof, let
$$
n_{k-1}<n-n^\alpha\leq n_k.
$$
Then
$$
n_{k-1}+2\left[\frac{n_{k-1}^\alpha}2\right]<n_{k-1}+n^\alpha<n,
$$
hence
$$
n_k+ \left[\frac{n_{k-1}^\alpha}2\right]<n.
$$
Furthermore  for $n$ large enough

$$\frac{n}{2}\leq n-n^{\alpha} \leq n_k$$
hence
$$\log n\leq \log (2 n_k) \leq\log ( 2n_{k+1})$$

 Thus with probability 1 for large $n$ the tube of radius
$\log n$ around the path $\{S_j,\, j=1,2,\ldots, n-[n^\alpha]\}$
is disjoint from the path $\{S_\ell,\, \ell=n,\ldots\}$, i.e.
Lemma 4.1 follows.

To prove Theorem 1.1 observe that it suffices to consider points visited
before time $n-n^\alpha$, ($2/(d-2)<\alpha<1$), since in the time
interval $(n-n^\alpha,n)$ the maximal local time is less than
$\alpha(1+\varepsilon)\lambda\log n$, hence this point cannot be in
${\cal A}_n$. Consequently, Theorem 1.1 follows from Theorem 1.2 and
Lemma 4.1.

\renewcommand{\thesection}{\arabic{section}.}

\section{Proof of Theorem 1.3}

\renewcommand{\thesection}{\arabic{section}} \setcounter{equation}{0}
\setcounter{theorem}{0} \setcounter{lemma}{0}

First we prove
\begin{lemma}
Let $A_i,B_i$ be events such that $\sum_i\bp(A_i)=\infty$,
$$
\bp(A_iA_k)\leq c_1\bp(A_i)\bp(A_k),
$$
and
$$
\bp(A_iB_i)\geq c_2\bp(A_i)
$$
with some constants $c_1,c_2>0$. Then
$$
\bp(A_iB_i\, {\rm i.o.})>0.
$$
\end{lemma}

\noindent{\bf Proof.}
$$
\sum_i\bp(A_iB_i)\geq c_2\sum_i\bp(A_i)=\infty.
$$
On the other hand,
$$
\bp(A_iB_iA_kB_k)\leq \bp(A_iA_k)\leq
\frac{c_1}{c_2^2}\bp(A_iB_i)\bp(A_kB_k),
$$
the Lemma follows by Borel-Cantelli lemma in Spitzer \cite{Sp}, pp. 317.

\bigskip
To prove the Theorem, define the stopping times $V_j$ as in R\'ev\'esz
\cite{R04}. Let
\begin{eqnarray*}
\rho_0(t)&=& t,\\
\rho_1(t)&=&\min\{\tau:\ \tau>t,\ S(\tau)=S(t)\},\\
\rho_2(t)&=&\min\{\tau:\ \tau>\rho_1(t),\ S(\tau)=S(\rho_1(t))=S(t)\},\\
\ldots, & &\\
\end{eqnarray*}
where here and the sequel we denote $S(k)=S_k$.
\begin{eqnarray*}
U(L,t)&=&\left\{\begin{array}{ll}
t+L & {\rm if\ } \rho_1(t)-t> L,\\
\rho_1(t)+L & {\rm if\ } \rho_1(t)-t\leq L,\ \rho_2(t)-\rho_1(t)> L,\\
\rho_2(t)+L & {\rm if\ } \rho_1(t)-t\leq L,\ \rho_2(t)-\rho_1(t)\leq L,\
\rho_3(t)-\rho_2(t)>L,\\
\ldots,&
\end{array}\right.\\
L_k&=&(\log(k+2))^\alpha,\quad (\alpha>\frac{2}{d-2},\ k=0,1,2,\ldots)\\
V_0&=&0, \quad V_{j+1}=U(L_{j},V_{j}), \quad (j=0,1,2,\ldots)
\end{eqnarray*}
$V_{j+1}$ is the first time-point after $V_j$ when the random walk has
not visited $S(V_j)$ during a time-interval of length $L_{j}$.

Let $\{x_n\}$ be a sequence of points in $Z^d$ as in Theorem 1.3
and define the events
\begin{equation}
A_j=\{\xi(S(V_j), V_{j+1})-\xi(S(V_j),V_j)=\psi_j,\,
\xi(S(V_j)+x_{V_j},V_{j+1})-\xi(S(V_j)+x_{V_j},V_j)=0\},
\label{aj}
\end{equation}
\begin{equation}
B_j=\{
\xi(S(V_j)+x_{V_j},V_j)=\xi(S(V_j)+x_{V_j},\infty)-
\xi(S(V_j)+x_{V_j},V_{j+1})=0\},
\label{bj}
\end{equation}
where $\psi_j=[\lambda(\log j+\log\log j)]$.

\begin{lemma}
The events $A_j$, $j=1,2,\ldots are$ independent and
\begin{equation}
\bp(A_j)\geq \frac{C_{14}}{j\log j}. \label{paj}
\end{equation}
\end{lemma}

\bigskip\noindent{\bf Proof.}
Since $\{V_j\}_{j=1}^\infty$ is a sequence of stopping times and $A_j$
depends only on the random walk between $V_j$ and $V_{j+1}$,
independence follows. To show (\ref{paj}), let $U_j:=U(L_j,0)$. Consider
the random walk starting from $V_j$ as a new origin. Then the original
random walk in the interval $(V_j,V_{j+1})$ has the same distribution as
the new random walk in $(0,U_j)$. Hence
$$
\bp(A_j\mid V_j=m)=\bp (\xi(0,U_j)=\psi_j,\, \xi(x_m,U_j)=0).
$$
The event $\{\xi(0,U_j)=\psi_j,\, \xi(x_m,U_j)=0\}$ means that there are
exactly $\psi_j$ excursions around $0$, each of which has length less
than $L_j$, none of them are visiting $x_m$ and in the last section
$(U_j-L_j,U_j)$ the random walk starting from $0$, does not visit $0$
and $x_m$. Hence applying (\ref{qxn}) of Lemma 2.2,
$$
\bp (\xi(0,U)=\psi_j,\, \xi(x_m,U)=0)
$$
$$
=\left(q_{x_m}+O((\log j)^{-\alpha(d/2-1)})\right)^{\psi_j}
\bp(\xi(0,L_j)=0,\, \xi(x_m,L_j)=0).
$$
Obviously
$$
\bp(\xi(0,L_j)=0,\, \xi(x_m,L_j)=0)\geq
\bp(\xi(0,\infty)=0,\, \xi(x_m,\infty)=0)= 1-q_{x_m}-s_{x_m}.
$$

From the inequalities (\ref{qe}) and (\ref{qse}) of Lemma 2.2 we can get
by easy calculation that
$$
\bp (\xi(0,U_j)=\psi_j,\, \xi(x_m,U_j)=0)\geq
C_{15}(q_{x_m})^{\psi_j} \geq C_{16}(1-\gamma)^{\psi_j}
\left(1-\frac{(1-\gamma_{x_m})^2}{1-\gamma}\right)^{\psi_j}.
$$
Since $L_j\geq 1$, we obviously have $V_j\geq j$, i.e. we can take
$m\geq j$. Since
$$
(1-\gamma)^{\psi_j}\geq \frac{1}{j\log j}
$$
and (cf. (\ref{greenas}))
$$
(1-\gamma_{x_m})^2\sim C_{17}(\log m)^{-1},
$$
we have
$$
\bp (A_j\mid V_j=m)=\bp (\xi(0,U_j)=\psi_j,\, \xi(x_m,U_j)=0)\geq
\frac{C_{14}}{j\log j},
$$
with $C_{14}>0$ independent of $m$, the lemma follows.

\begin{lemma} Let the events $A_j$, $B_j$ be defined by {\rm (\ref{aj})}
and {\rm (\ref{bj})}. Then
\begin{equation}
\bp(A_jB_j)\geq \gamma^2\bp(A_j).
\label{ajbj}
\end{equation}
\end{lemma}

\bigskip\noindent{\bf Proof.}
$$
\bp(A_jB_j)=\be \bp(A_jB_j\mid S(V_j),\, S(V_{j+1}))
$$
$$
=\be\left(\bp(A_j\mid S(V_j),\, S(V_{j+1}))\bp(B_j\mid S(V_j),\,
S(V_{j+1}))\right).
$$
We show that
\begin{equation}
\bp(B_j\mid S(V_j),\, S(V_{j+1}))\geq \gamma^2,\quad j=1,2,\ldots
\label{pbj}
\end{equation}
Consider the reversed random walk before $S(V_j)$, as in the the
proof of Theorem 1.2, i.e. $S_i'=S(V_j-i)-S(V_j)$, and its local time
$\xi'(x,n)$ and also the forward random walk starting from $S(V_{j+1})$,
i.e. $S_i"=S(V_{j+1}+i)-S(V_{j+1})$, $i=1,2,\ldots$ and its local time
$\xi"(x,n)$. These two random walks are independent and the event $B_j$
means that the first random walk $S'$ does not visit $x_{V_j}$ (up to
time $V_j$) and the second random walk $S"$ does not visit
$S(V_j)+x_{V_j}-S(V_{j+1})$ (for infinite time). Hence
\begin{eqnarray*}
&&\bp(B_j\mid S(V_j),\, S(V_{j+1})) \\
\quad &=&\bp(\xi'(x_{V_j},V_j)=0,\,
\xi"(S(V_j)-S(V_{j+1})+x_{V_j},\infty)=0\mid
S(V_j),S(V_{j+1})) \\
\quad &\geq&\bp(\xi'(x_{V_j},\infty)=0)
\bp(\xi"(S(V_j)-S(V_{j+1})+x_{V_j},\infty)=0\mid
S(V_j),S(V_{j+1})).
\end{eqnarray*}
From (\ref{xe}) of Lemma 2.2 it follows that
$$
\bp(\xi'(x_{V_j},\infty)=0)\geq \gamma
$$
and similarly
$$
\bp(\xi"(S(V_j)-S(V_{j+1})+x_{V_j},\infty)=0\mid
S(V_j),S(V_{j+1}))\geq \gamma,
$$
hence (\ref{pbj}) follows, which, in turn, implies (\ref{ajbj}). This
proves Lemma 5.3.

Lemma 5.2 and Lemma 5.3 together imply by Lemma 5.1 that
$$
\bp(A_jB_j\, {\rm i.o.})>0.
$$
Since (cf. R\'ev\'esz \cite{R04})
$$
V_{j}=n_j\leq O(1)j(\log j)^\alpha\quad {\rm a.s.},
$$
assuming that $A_jB_j$ occurs, we have
\begin{eqnarray*}
\xi(S_{n_j},\infty)&=&
\xi(S(V_{j+1}),\infty)\geq\xi(S(V_j),V_{j+1})-\xi(S(V_j),V_j)
\geq\psi_j\geq\\
&\geq&\lambda\log n_j-\lambda\alpha\log\log n_j+(1-\varepsilon)
\lambda\log\log n_j\geq\\
&\geq&\lambda\log n_j+\lambda\left({\frac{d-4}{d-2}}-\varepsilon\right)
\log\log n_j
\end{eqnarray*}
and also $\xi(S_{n_j}+x_{n_j},\infty)=0$. Thus we have $\bp(D_n\, {\rm
i.o.})>0$, where
$$
D_n=\left\{\xi(S_n,\infty)\geq \lambda\left(\log n
+\left(\frac{d-4}{d-2}-\varepsilon\right)\log\log n\right),\quad
\xi(S_n+x_n,\infty)=0\right\}.
$$
Let
\begin{eqnarray*}
\widetilde D_n=&&\Big\{\xi(S_n,\infty)\geq \lambda\left(\log n
+\left(\frac{d-4}{d-2}-\varepsilon\right)\log\log n\right),\\
&&\xi(S_n+x_n,\infty)-\xi(S_n+x_n,\log n)=0\Big\}.
\end{eqnarray*}
Then we have also $\bp(\widetilde D_n\, {\rm i.o.}) >0$ and since
$\widetilde D_n$ is a tail event for the random walk, by 0-1 law we have
$\bp(\widetilde D_n\, {\rm i.o.})=1$.

To show that also $\bp(D_n\, {\rm i.o.})=1$, we prove the following
\begin{lemma}
For any $0<\delta<1/2$ with probability 1 there exists $n_0$ such that
for $n\geq n_0$ we have
$$
\xi(S_n+x,n^{\delta})=0\quad \text{\rm for all}\quad |x|\leq \log n.
$$
\end{lemma}
{\bf Proof.} By (\ref{jp}) we get
$$
\bp\left(\bigcup_{|x|\leq \log n}\bigcup_{j\leq
n^{\delta}}\{S_j=S_n+x\}\right) \leq \sum_{|x|\leq\log
n}\sum_{j\leq n^{\delta}}\bp(S_j=S_n+x)
$$
$$
\leq \sum_{|x|\leq \log n}\sum_{j\leq n^{\delta}}\frac{C_8}{(n-j)^{d/2}}
\leq \frac{C_{17}(\log n)^{d}}{n^{d/2-\delta}},
$$
and since this is summable, the lemma follows by Borel-Cantelli lemma.
This implies $\bp(D_n\, {\rm i.o.})=1$, proving Theorem 1.3.

\end{document}